\documentclass[12pt]{amsart}
\usepackage{amscd,amssymb,verbatim}
\usepackage{amssymb,amsmath,amsthm,epsfig}

\def\cnum#1{\bigcirc\kern -8pt#1}
\newcommand{\n}{\noindent}
\newcommand{\vp}{\varepsilon}

\theoremstyle{plain}
\newtheorem{thm}{Theorem}[section]
\newtheorem{fact}[thm]{Fact}
\newtheorem{cor}[thm]{Corollary}

\newtheorem{lem}[thm]{Lemma}
\newtheorem{defn}[thm]{Definition}
\newtheorem{rem}[thm]{Remark}

\newcommand{\T}{\mathcal T}

\newcommand{\R}{\mathcal R}
\newcommand{\NC}{\mathcal N}
\newcommand{\supp}{{\rm supp}\,}
\newcommand{\spn}{{\rm span}\,}

\newcommand{\sign}{{\rm sign}\,}

\newcommand{\Ker}{{\rm Ker}\,}
\newcommand{\N}{\mathbb{N}}

\newcommand{\Q}{\mathbb{Q}}

\newcommand{\mnorm}[1]{\left\vert\kern-0.9pt\left\vert\kern-0.9pt\left\vert #1
    \right\vert\kern-0.9pt\right\vert\kern-0.9pt\right\vert}

\begin{document}

\title[Operators on Asymptotic $\ell_p$ Spaces]{Operators on Asymptotic $\ell_p$ Spaces which are not Compact Perturbations of a Multiple of the Identity}

\author{Kevin Beanland}
\date{}
\thanks{The author would like to thank G. Androulakis for his guidance
during the preparation of this manuscript.}

\begin{abstract}
We give sufficient conditions on an asymptotic $\ell_p$ (for $1 < p
< \infty$) Banach space which ensure the space admits an operator
which is not a compact perturbation of a multiple of the identity.
These conditions imply the existence of strictly singular
non-compact operators on the HI spaces constructed by G. Androulakis
and the author and by I. Deliyanni and A. Manoussakis. Additionally
we show that under these same conditions on the space $X$,
$\ell_\infty$ embeds isomorphically into the space of bounded linear
operators on $X$.
\end{abstract}

\maketitle

\markboth{K. Beanland}{Operators on Asymptotic $\ell_p$ Spaces which
are not Compact Perturbations of a Multiple of the Identity}

\section{Introduction} \label{sec1}

In this note, we give sufficient conditions on
a Banach space whereby the space of bounded linear operators does does not solve the
scalar-plus-compact problem of Lindenstrauss.  

Lindenstrauss' question is related to the result of N. Aronszajn and
K.T. Smith (\cite{AS}) in 1954, which implies that if a space $X$
satisfies the above condition and is a complex space, then every
bounded linear operator on $X$ must have a non-trivial invariant
subspace.  Thus, a complex space which is a positive solution to
Lindenstrauss' problem also serves as a positive solution to the
invariant subspace problem for Banach spaces.

That being said, the possibility that for any Banach space there is
an operator on the space which is not a compact perturbation of a
multiple of the identity, is still in play.  In support of this
possibility, sufficient conditions have been established on a space
$X$ which imply $\ell_\infty$ embeds isomorphically into
$\mathcal{L}(X)$, the space of bounded linear operators on $X$ (see
\cite{ABDS},\cite{Fe},\cite{K}).  If a space $X$ serves as a
positive solution to the scalar plus compact problem, has a basis
(or more generally the Approximation Property) and a separable dual
space, then $\mathcal{L}(X)$ must be separable. Curiously, each of
the results in support of a negative solution to the scalar plus
compact problem, require the existence of an unconditional basic
sequence in the space.  The weaker problem of whether there is an
operator which is not a compact perturbation of a multiple of the
inclusion from a subspace of a Banach space to the whole space has
also received attention (see
\cite{AOST},\cite{ASa},\cite{Go1},\cite{S}).

In their successful effort to construct the first example of a space
with no unconditional basic sequence, W.T. Gowers and B. Maurey
(\cite{GM1}) constructed a space which, as W.B. Johnson observed,
possesses a stronger property called hereditarily indecomposable
(HI). A Banach space is HI if no (closed) infinite dimensional
subspace can be decomposed into a direct sum of two further infinite
dimensional subspaces.  This groundbreaking construction was a great
leap forward in the progression towards a positive solution to the
scalar-plus-compact problem. More precisely, it was shown that every
operator on the space of Gowers-Maurey can be decomposed as a
strictly singular perturbation of a multiple of the identity
operator. Spaces which have this property are now aptly referred to
as spaces admitting ``few operators.''   An operator on a Banach
space is called strictly singular if the restriction of it to any
infinite dimensional subspace is not an isomorphism. The ideal of
strictly singular operators on a space contains that of the compact
operators, but in some cases (e.g. $\ell_p$, $1 \leq p < \infty$)
they coincide. The fact that Gowers-Maurey space admits few
operators is related to the fact that it is HI. In fact it was shown
in \cite{GM1}, that every complex HI space admits few operators.  In
1997, V. Ferenczi proved (\cite{F1}) that a complex space $X$ is HI
if and only if every operator from a subspace of $X$ into $X$ is a
multiple of the inclusion plus a strictly singular operator.  It is
not the case however, that admitting few operators implies that the
space is HI. The most recent in a collection of counterexamples is
the paper of S.A. Argyros and A. Manoussakis (\cite{AM}) in which
they construct a reflexive space admitting few operators for which
every Schauder basic sequence has an unconditional subsequence. The
most comprehensive resource for HI spaces and spaces admitting few
operators is \cite{AT}.

The natural question then becomes: for any of these spaces which
admit few operators does there exist a strictly singular non-compact
operator, or do the strictly singular and compact ideals coincide?
There have been results in this direction as well.  In 2000, Argyros
and Felouzis (\cite{AF}) constructed an HI space $X$ with the
property that for every infinite dimensional subspace of $X$ there
is a strictly singular non-compact operator on $X$ with range
contained in the subspace. In 2001, G. Androulakis and Th.
Schlumprecht (\cite{AS}) constructed a strictly singular non-compact
operator on the space of Gowers-Maurey. In 2002, I. Gasparis
(\cite{Ga1}) did the same for certain members of the class of
totally incomparable asymptotic $\ell_1$ HI spaces constructed in
\cite{Ga2}.  In 2006, G. Androulakis and the author (\cite{AB}) and
A. Manoussakis and I. Deliyanni (\cite{DM}) independently
constructed different asymptotic $\ell_p$ HI spaces (for $p=2$ in
the former case and  for all $1 < p < \infty$ in the latter). In the
following, by extending results in \cite{Ga1}, sufficient conditions
are established under which a strictly singular non-compact
operators can be found on each of these spaces.\\

\n {\bf Additional Note:} Recently, S.A. Argyros and R. Haydon have constructed 
a $\mathcal{L}_\infty$ HI space on with every operator is a multiple of the identity plus a compact operator.

\section{Definitions and Notation} \label{sec2}

Our notation is standard and can be found in \cite{LT}.  Let
$(e_i)_{i=1}^\infty$ denote the unit vector basis of
$c_{00}(\N)=c_{00}$, and $(e^*_i)_{i=1}^\infty$ the biorthogonal
functionals of $(e_i)_i$. Let $\spn\{(e_i)_i\}$ denote vectors
finitely supported on $(e_i)_i$. For a Banach space $X$ let
$Ba(X)=\{x\in X : \|x\| \leq 1 \}$ and $S(X)=\{x \in X : \| x\| =1
\}$. If $E,F \subset \N$ then $E<F$ if $\max E < \min F$. If $x =
\sum^\infty_{i=1} a_i e_i$ for scalars $(a_i)_i$, let $\supp(x)= \{i
: a_i \not=0\}$ and the range of $x$, denoted $r(x)$, be the
smallest interval containing $\supp(x)$.

The notion of Schreier families (\cite{AA}) is used throughout. They
are defined inductively as follows.  Let $S_0=\{\{n\}: n \in
\N\}\cup \{\emptyset\}$.  After defining $S_n$ let,

$$S_{n+1}=\{F \subset \N : F= \bigcup_{i = 1}^m F_i ~\mbox{for}~
F_i \in S_n ~\mbox{and}~ m \in \N, ~ m \leq F_1 < \ldots < F_m \}
\cup \{\emptyset\}$$

\n A few properties of the Schreier families we need are:
\begin{itemize}
\item[$\cdot$](Hereditary) For $n \in \N$, $S_n \subset S_{n+1}$.
\item[$\cdot$](Spreading) If $(p_i)_{i=1}^N \in S_n$ and $p_i \leq
q_i$ for all $i \leq N$, then $(q_i)_{i=1}^N \in S_n$.
\item[$\cdot$](Convolution) If $(F_i)_{i=1}^N$ is a collection of subsets of $\N$ such that $F_i
\in S_n$ for all $i\leq N$, $F_1 < \ldots < F_N$ and $(\min
F_i)_{i=1}^N \in S_m$ for some $n,m \in \N$ then $\bigcup_{i=1}^N
F_i \in S_{n+m}$.
\end{itemize}

\n Let $(E_i)_{i=1}^k$ be a sequence of successive subsets of $\N$,
we say that $(E_i)_{i=1}^k$ is $S_n$ admissible if $(\min
E_i)_{i=1}^k \in S_n$.  For $E_1 < \cdots < E_k \subset \N$ and
$(a_j)_j \in c_{00}$ the sequence $(x_i)_{i=1}^k$ defined by $x_i=
\sum_{j\in E_i} a_j e_j$ is called a block sequence of $(e_j)_j$.
For a block sequence $(x_i)_{i=1}^k$ of $(e_j)_j$ we say that
$(x_i)_{i=i}^k$ is $S_n$ admissible if $(\supp x_i)_{i=1}^k$ is
$S_n$ admissible.

Herein we define a class of spaces in terms of the norming
functionals of the space. We begin by recalling the notion of a
norming set (\cite{Ga1}),

\begin{defn}
A set $\NC \subset \spn\{(e_i^*)_i\}$ is called norming if the
following conditions hold,
\begin{itemize}
\item[$\cdot$] $(e^*_n)_n \subset \NC$,
\item[$\cdot$] If $x^* \in \NC$ then $|x^*(e_n)|\leq 1$ for all $n \in
\N$.
\item[$\cdot$] If $x^* \in \NC$ then $-x^* \in \NC$ ($\NC$ is symmetric).
\item[$\cdot$] If $x^*  \in \NC$ and $E$ is an interval in $\N$ then $Ex^* \in
\NC$  (where $Ex^*$ denotes the restriction of $x^*$ to the
coordinates in $E$).
\end{itemize}
\end{defn}

\n If $\NC$ is a norming set we can define a norm $\|\cdot \|_{\NC}$
on $c_{00}$ by

$$\biggl\|\sum_i a_i e_i \biggr\|_{\NC} = \sup \biggl\{ x^*\biggl(\sum_i a_i e_i\biggr) : x^* \in
\NC \biggr\}$$

\n for every $(a_i) \in c_{00}$.  Now define the Banach space
$X_{\NC}$ to be the completion of $c_{00}$ under the above norm.  By
the definition of norming set, $(e_i)_i$ is a normalized bimonotone
basis for $X_{\NC}$.

For the following definitions and notation we closely follow
\cite{Ga1}.  The following are conditions on two increasing
increasing sequences of positive integers, $(n_i)_{i=1}^\infty$ and
$(m_i)_{i=1}^\infty$.

\begin{itemize}
\item[(i)] $m_1 > 3$, there is an increasing sequence of
positive integers $(s_i)_{i=1}^\infty$ such that
$m_{2j}=\prod_{i=1}^{j-1} m_{2i}^{s_i}$, $m_{2j+1}=m_{2i}^5$ for $i
\geq 1$ and $m_1^5=m_2$.

\item[(ii)] For the sequence of integers, $(f_i)_{i=2}^\infty$ defined
by,

$$f_j = \max
\left\{ \rho n_1 +\sum_{1 \leq i<j}\rho_i n_{2i} : \rho, \rho_i \in
\mathbb{N}\cup \{0\}, m_1^\rho \prod_{1 \leq i<j}m_{2i}^{\rho_i} <
m_{2j} \right\}$$

\n require that $4f_j < n_{2j}$ for all $j \geq 2$ and $5n_1 < n_2$.
\end{itemize}

\n We now define a particular type of norming set.  Our definition
is slightly less general than that which would be considered
analogous to $(M,N)$-Schreier in \cite{Ga1}. Our goal is to tailor
the definition of $(M,N,q)$-Schreier so as to make it as apparent as
possible that the spaces found in \cite{AB} and \cite{DM} are
$(M,N,q)$-Schreier for specified $q$.

\begin{defn}
For sequences $M=(m_i)_{i=1}^\infty$ and $N=(n_i)_{i=1}^\infty$
satisfying (i) and (ii) we call a norming set $\NC$,
$(M,N,q)$-Schreier (for $1/q+1/p=1$ and $1 < p,q < \infty$) if for
the following sets, with $k \in \N$,

$$ \NC_k = \biggl\{\frac{1}{m_{2k}} \sum_i \gamma_i x^*_i : (\gamma_i)_i \in
Ba(\ell_q), ~ \gamma_i \in \Q, ~(x^*_i)_i ~\mbox{is}~ S_{n_{2k}}
\mbox{admissible and} ~ (x^*_i)_i \subset \NC \biggr\} $$

\begin{equation*}
\begin{split}
\NC_{\infty}^q = \bigcup_{k=0}^\infty \biggl\{ \frac{1}{m_{2k+1}}
\sum_i \gamma_i E x^*_i : & ~(\gamma_i)_i \in 2^{1/p} Ba(\ell_q),
~\gamma_i \in \Q,~ E~\mbox{is an interval $\subset \N$,} ~ (x_i^*)_i
~
 \\
& \mbox{is an}~ S_{n_{2k+1}}\mbox{admissible} \subset
\bigcup_{j=1}^\infty \NC_{j} \}
\end{split}
\end{equation*}

\n we have $\NC \subset \bigcup_{j=1}^\infty \NC_j \cup
\NC_{\infty}^q \cup \{\pm e_n : n \in \N \}$ and $\NC_j \subset \NC$
for all $j$.

\end{defn}

In the case of the Banach space constructed in \cite{AB}, there are
fixed sequences $M=(m_i)_{i=1}^\infty$ and $N=(n_i)_{i=1}^\infty$.
We suppose further that $M$ and $N$ satisfy conditions (i), (ii). It
follows directly from the definition that the norming set for this
space is $(M,N,2)$-Schreier.

For the asymptotic $\ell_p$ HI space, $X_{(p)}$, found in \cite{DM}
the reasoning is similar.  Assume that the sequences $M$ and $N$
prescribed in \cite{DM} satisfy conditions (i), (ii).  For a fixed
$p$ and $1/q+1/p=1$ we must show that the norming set $\NC$ (denoted
$K$ in \cite{DM}) is $(M,N,q)$-Schreier.  The reader should refer to
\cite{DM} for the precise definitions of $K,K^n$ and $K_j^n$.  $K=
\bigcup_{n=1}^\infty K^n$ where $K^n =\bigcup_{j=1}^\infty K_j^n$.
By definition, for $j,n \in \N$, $K_{2j}^n \subset \NC_j$, and
$K_{2j+1}^n \subset \NC_\infty^q$ (for the latter inclusion the
factor $2^{1/p}$ in the definition of $\NC_\infty^q$ is required).
Thus, $K \subset \bigcup_j \NC_{j}\cup \NC_{\infty}^q \cup \{\pm e_n
: n \in \N \}$. If $x^* \in \NC_j$ then $x^* = 1/m_{2j} \sum_i
\gamma_i x^*_i$ where $(x^*_i)_i \subset \NC = K$, so again by the
definition of $K$, we see that $\NC_j \subset \NC$.

It is convenient to view an element of $\NC$ as successive blocks of
the basis $(e_i^*)_i$. This decomposition into blocks is not unique,
and thus our goal is to find a decomposition that is the most
suitable. To this end, we associate each element of $\NC$ with a
rooted tree. A finite set with a partial ordering $(\T,\preceq)$ is
called a tree if for every $\alpha \in \T$ the set $\{ \beta \in \T
: \beta \preceq \alpha \}$ is linearly ordered. Each element of the
tree $\T$ is called a node. A node $\alpha \in \T$ such that there
is no $\beta$ with $\alpha \prec \beta$ is called terminal ($\alpha
\prec \beta$ means $\alpha \preceq \beta$ and $\alpha \not= \beta$).
If $\beta \prec \alpha$ we say $\alpha$ is a successor of $\beta$.
For $\alpha \in \T$ let $D_{\alpha}(\T)$ denote the set of immediate
successors of $\alpha$ in $\T$.  A branch of $\T$ is a maximal
linearly ordered subset.

For each $\alpha \in \T$ we define corresponding $\gamma_\alpha \in
\Q$, $m_\alpha \in (m_i)_{i=1}^\infty$ and $n_\alpha \in
(n_i)_{i=1}^\infty$.  We associate to each $x^* \in \NC$ a rooted
tree $\T$ (i.e. a tree with a unique first node) in the following
way:  Let $\alpha_0$ be the root of $\T$. There is an
$x^*_{\alpha_0} \in \NC$ such that $x^*=\gamma_{\alpha_0}
x^*_{\alpha_0}$ and

$$ x^*_{\alpha_0} = \frac{1}{m_{\alpha_0}} \sum_{\beta \in
D_{\alpha_0}(\T)} \gamma_\beta x^*_\beta. $$

\n In this definition $(x^*_\beta)_{\beta \in D_{\alpha_0}(\T)}$ is
$S_{n_{\alpha_0}}$ admissible, $m_{\alpha_0}= m_j$ for some $j \in
\N$ and $(\gamma_\beta)_{\beta \in D_{\alpha_0}(\T)} \in Ba(\ell_q)$
if $j$ is even and $(\gamma_\beta)_{\beta \in D_{\alpha_0}(\T)} \in
2^{1/p}Ba(\ell_q)$ if $j$ is odd.  Thus for any pairwise
incomparable collection $A$ of $\T$ which intersects every branch of
$\T$ we have,

\begin{equation}
x^* = \sum_{\alpha \in A} \frac{\prod_{\beta \preceq \alpha}
\gamma_\beta}{\prod_{\beta \prec \alpha} m_\beta} x^*_{\alpha}.
\label{decomp}
\end{equation}

\n Call $(x^*_\alpha)_{\alpha \in \T}$ the functional tree of $x^*$.
For $\beta \in \T$ the functional $x_\beta^*$ has a corresponding
tree $\T_\beta$, which is a subset of $\T$.

\section{Main Results} \label{sec3}

We start this section by stating the main theorem of the paper. The
proof of this theorem can be found at the end of this section.  The
majority of this section is devoted to proving auxiliary lemmas and
remarks.

\begin{thm}
Let $\NC$ be a norming set which is $(M,N,q)$-Schreier and $X_\NC$
be the corresponding Banach space. There is an operator on $X_\NC$
which is not a compact perturbation of a multiple of the identity.
If $X_\NC$ is HI, this operator is strictly singular. Moreover,
$\ell_\infty$ embeds isomorphically into $\mathcal{L}(X_\NC)$.
\label{mainthm}
\end{thm}

The existence of the following sequence in the space $X_\NC^*$ is
the main ingredient in the construction of the desired operator. The
definition below is tailored to fit our construction.

\begin{defn}
Let $(x_k)_k$ be a block basic sequence of $(e_k)_k$.  If there is a
$C>0$ such that for all $l \in \N$, $F \subset \N$ with $F \geq l$
and $(x_k)_{k \in F}$ being $S_{f_l}$ admissible we have $\|\sum_{k
\in F} \beta_k x_k \| \leq C \|(\beta_{k})_{k \in F}\|_q$ for every
scalar sequence $(\beta_{k})_{k \in F}$, we say $(x_k)_k$ satisfies
an upper $\ell_q^\omega$ estimate with constant $C$.
\end{defn}

For the rest of the section we fix $p$, $q$ and $\NC$ such that $1/p
+ 1/q =1$ and $\NC$ is a $(M,N,q)$-Schreier norming set. It follows
easily that any normalized block sequence $(x_i)_{i=1}^m$ with $m
\leq \supp x_1$ in $X_\NC$, satisfies a lower $\ell_p$ estimate with
constant $1/m_2$. The next remark demonstrates that $X_\NC$ is an
asymptotic $\ell_p$ space by verifying that is satisfies an upper
$\ell_p$ estimate on normalized blocks.

\begin{rem}
Let $(x_i)_{i=1}^m \in X_\NC$ be a normalized block basic sequence
of $(e_i)_i$.  For any sequence of scalars $(a_i)_i$ the following
holds:

$$\biggl\|\sum_{i=1}^m a_i x_i \biggr\|_{\NC} \leq 12 \biggl( \sum_{i=1}^m |a_i|^p\biggr)^{\frac{1}{p}}. $$
\label{ulp}
\end{rem}

\begin{proof}

For $x^* \in \NC$, let $o(x^*)$ denote the height (i.e. the length
of the longest branch) of the tree $\T$ associated with $x^*$.  We
proceed by induction on $o(x^*)$. We will show that for all $x^* \in
\NC$ such that $o(x^*)=n$, $x^*(\sum_{i=1}^m a_i x_i) \leq 12 (
\sum_{i=1}^m |a_i|^p)^{\frac{1}{p}}$ holds for any normalized block
basic sequence $(x_i)_{i=1}^m$ of $(e_i)_i$ and any sequence of
scalars $(a_i)_i$. For $x^* \in \NC$ such that $o(\T)=1$ the
assertion follows easily. Assume the claim for all $y^* \in \NC$
such that $o(y^*)<n$ and let $x^*=1/m_k \sum_j \gamma_j x^*_j \in
\NC$ with $o(x^*)=n$. By definition of $\NC$, $(x^*_j)_j$ is
$S_{n_k}$ admissible and $(\gamma_j)_j \in 2^{1/p}Ba(\ell_q)$.
Define the following two sets,

$$Q(1)=\{1\leq i \leq m:~\mbox{there is exactly one}~ j ~\mbox{such that}~ r(x_j^*)\cap r(x_i) \not= \emptyset \},$$

\n and $Q(2)=\{1,\ldots , m\}\setminus Q(1)$.  Apply the functional
$x^*$ to $\sum_{i=1}^m a_i x_i$ to obtain:

\begin{equation*}
\begin{split}
\biggl| \frac{1}{m_k}\sum_{j} \gamma_j x^*_j \biggl(\sum_{i=1}^m a_i
x_i \biggr) \biggr|  & \leq \frac{1}{m_k}\sum_{j} |\gamma_j| \biggl|
x^*_j \sum_{\substack{i\in Q(1)\\ r(x_j^*)\cap r(x_i) \not=
\emptyset}} a_i x_i \biggr|  +
\frac{1}{m_k} \biggr| \sum_{j} \gamma_j x^*_j \sum_{i \in Q(2)} a_i x_i \biggr|  \\
    & \leq \frac{12}{m_k}\sum_j | \gamma_j |  \biggl
    (\sum_{\substack{i\in Q(1) \\ r(x_j^*)\cap r(x_i) \not= \emptyset}} |a_i|^p\biggr)^{\frac{1}{p}}
    + \sum_{i \in Q(2)} |a_i|\biggl| \frac{1}{m_k}\sum_j \gamma_j  x^*_j(x_i)\biggr| ,
    \label{eqn:CC}
\end{split}
\end{equation*}

\n The first inequality follows from the triangle inequality. The
second follows from applying the induction hypothesis for
$x_j^*(\sum_{\{i \in Q(1): r(x^*_j) \cap r(x_i) \not= \emptyset \}}
a_i x_i)$ and using the definition of $Q(2)$.  We may apply the
induction hypothesis since the height of the trees associated with
the functionals $x_j^*$ are each less than $n$. Before continuing,
notice that for each $i \in Q(2)$ the set $J_i=\{j : r(x_j^*)\cap
(x_i) \not= \emptyset\}$ is an interval and therefore,

\begin{equation}
\frac{1}{m_k} \sum_{j \in J_i} \frac{\gamma_j}{(\sum_{j\in J_i}
|\gamma_j|^q)^{1/q}} x^*_j \in \NC \label{inNC}
\end{equation}

\n The above estimate continues as follows,

\begin{equation*}
\begin{split}
&\leq 4 \sum_j |\gamma_j| \biggl(\sum_{\substack{i\in Q(1) \\
r(x_j^*)\cap r(x_i) \not= \emptyset}} |a_i|^p \biggr)^{\frac{1}{p}}
+
\sum_{i \in Q(2)} |a_i| \biggl( \sum_{j \in J_i} |\gamma_j|^q \biggr)^{\frac{1}{q}} \\
    & \leq 4 \left(\sum_j|\gamma_j|^q\right)^\frac{1}{q} \biggl(\sum_j \sum_{\substack{i\in Q(1) \\
    r(x_j^*)\cap r(x_i) \not= \emptyset}} |a_i|^p\biggr)^{\frac{1}{p}} +
    \biggr( \sum_{i \in Q(2)} |a_i|^p \biggr)^\frac{1}{p} \left(\sum_{i\in Q(2)} \sum_{j \in J_i}
    |\gamma_j|^q\right)^{\frac{1}{q}} \\
    & \leq 4 \left(\sum_j|\gamma_j|^q\right)^\frac{1}{q} \biggl(\sum_j \sum_{\substack{i\in Q(1) \\
    r(x_j^*)\cap r(x_i) \not= \emptyset}} |a_i|^p\biggr)^{\frac{1}{p}}
    + \biggl( \sum_{i \in Q(2)} |a_i|^p \biggr)^\frac{1}{p}\biggl(2 \sum_{j} |\gamma_j|^q\biggr)^\frac{1}{q}  \leq
    12 \biggl(\sum_{i}^m |a_i|^p\biggr)^{\frac{1}{p}}.
    \label{eqn:BB}
\end{split}
\end{equation*}

\n In the first inequality we used the fact that $3<m_1$ in the
first term and (\ref{inNC}) in the second term. For the second
inequality we applied H\"{o}lders inequality. For the third
inequality we used the fact that for each $j$ there are at most two
values of $i \in Q(2)$ such that $r(x^*_j) \cap r(x_i) \not=
\emptyset$. For the final inequality we used $(\sum_\ell
|\gamma_\ell|^q)^{1/q} \leq 2$. This finishes the proof.
\end{proof}

The following is a compilation of remarks (variants of which can be
found in \cite{Ga1}) regarding the sequences
$(m_i)_{i=1}^\infty,(n_i)_{i=1}^\infty$ and $(f_i)_{i=2}^\infty$. In
the interest of completeness we have included the proofs.
\begin{itemize}
\item[(1.1)] If $p_k =5n_1+\sum_{i<k} s_i n_{2i}$ for $k \geq 2$, then $p_k \leq 2
f_k$.
\item[(1.2)] If $(a_i)_{i=1}^{k-1}$ is a sequence of non-negative integers
and $a \in \N \cup \{0\}$ such that $m_1^a \prod_{i<k} m_{2i}^{a_i}
< m_{2k}$ then $an_1+ \sum_{i<k} a_i n_{2i} < p_k$.
\item[(1.3)] Let $(a_\ell)_{\ell=1}^{k-1}$ be a sequence of non-negative integers,
$(x_i^*)_{i=1}^t \in \NC$ be $S_{\sum_{l<k} a_l n_{2l}}$ admissible
and $(\beta_i)_{i=1}^t \in Ba(\ell_q)$ then we have,

$$\frac{1}{\prod_{\ell<k} m_{2\ell}^{a_\ell}} \sum_{i=1}^t \beta_i x^*_i \in \NC .$$
\end{itemize}

The proof of (1.1) follows by induction.  For $k=2$ we have
$f_2=4n_1+(s_1-1)n_2$.  Since $s_1 \geq 2$ the claim follows.
Suppose the statement is true for some $k \geq 2$. Let $f_{k}=\gamma
n_1+ \sum_{i< k} \gamma_i n_{2i}$ and observe that,

\begin{equation*}
\begin{split}
p_{k+1}  = 5n_1 + \sum_{i < k} s_i n_{2i} +s_{k}n_{2k}
 & \leq 2 (\gamma n_1 + \sum_{i < k} \gamma_i n_{2i}) + s_{k} n_{2k} ~(\mbox{by
 the induction hypothesis}) \\
 & \leq 2 (\gamma n_1 + \sum_{i < k} \gamma_i n_{2i} + s_{k} n_{2k}) \leq 2
 f_{k+1}.
\end{split}
\end{equation*}

\n We obtained the third inequality by noting that $m_1^\gamma
\prod_{i<k}m_{2i}^{\gamma_i}m_{2k}^{s_{k}} < m_{2k} m_{2k}^{s_k} =
m_{2k+2}$ and using the maximality of $f_{k+1}$.

To prove (1.2), again proceed by induction.  For $k=2$, deduce from
the hypothesis that $a+5a_1 < 5s_1$.  Clearly, $a_1 < s_1$.  If
$a<5$ we are done.  Suppose $5n \leq a < 5(n+1)$ for some $n \in
\N$. This implies that $a_1 < s_1-n$.  The following inequality
finishes the proof of the base step,

$$a n_1 + a_1 n_2 < a n_1 + (s_1 - n)n_2 \leq s_1 n_2 +5(n+1)n_1 -n
n_2 < 5 n_1 + s_1 n_2. $$

\n The final inequality follows from $5n_1 < n_2$.

Assume the statement is true for some $k\geq 2$. By assumption,
$m_1^a\prod_{i<k+1}m_{2i}^{a_i} < m_{2k+2}$ and by definition
$m_{2k+2} = m_{2k}^{s_k+1}$. Clearly, $m_1^a\prod_{i<k}m_{2i}^{a_i}
< m_{2k}^{s_k-a_k+1}$. Thus $s_k+1 \geq a_k$. This leaves two
possibilities, either $s_k = a_k$ or $s_k > a_k$. In the former
case, $m_1^a \prod_{i<k}m_{2i}^{a_i} < m_{2k}$. By the induction
hypothesis $an_1 + \sum_{i<k} a_i n_{2i} < p_{k}$ and thus $an_1 +
\sum_{i<k+1} a_i n_{2i} < p_{k+1}$. If $s_k
> a_k$ we claim that $an_1 + \sum_{i<k+1} a_i n_{2i} < s_k n_{2k}$, which
clearly finishes the proof. To see this, we start by showing that $a
n_1 + \sum_{i<k} a_i n_{2i} \leq 2(s_k-a_k+1)f_k$.  By assumption
$m_1^a \prod_{i<k}m_{2i}^{a_i} < m_{2k}^{s_k-a_k+1}$, which implies
that,

$$\displaystyle m_1^{\left\lfloor\frac{a}{s_k-a_k+1}\right\rfloor} \prod_{i<k}m_{2i}^{\left\lfloor\frac{a_i}{s_k-a_k+1}\right\rfloor} < m_{2k}.$$

\n where $\lfloor x \rfloor$ is the greatest integer of $x$.  By the
maximality of $f_k$ we have,

$$\left\lfloor\frac{a}{s_k-a_k+1}\right\rfloor n_1 + \sum_{i<k} \left\lfloor\frac{a_i}{s_k-a_k+1}\right\rfloor n_{2i} \leq f_k. $$

\n Since $x \leq 2 \lfloor x \rfloor$ for $x \geq 0$ we see that,

$$\frac{a}{s_k-a_k+1} n_1 + \sum_{i<k} \frac{a_i}{s_k-a_k+1} n_{2i} \leq 2
\biggl(\left\lfloor\frac{a}{s_k-a_k+1}\right\rfloor n_1 +
\sum_{i<k} \left\lfloor\frac{a_i}{s_k-a_k+1}\right\rfloor
n_{2i}\biggr) \leq 2f_k.$$

\n Finally, using $4f_k < n_{2k}$ observe that,

$$a n_1 + \sum_{i<k+1} a_i n_{2i} \leq 2(s_k-a_k+1)f_k +a_k n_{2k} <
n_{2k}((s_k+1)/2-a_k) +a_k n_{2k} < s_k n_{2k}.$$

The proof of (1.3) requires a complicated induction.  For simplicity
we prove the case where $a_j=a_l=1$ for some $j,l \leq k$. Suppose
$(x^*_i)_{i=1}^t \in \NC$ is $S_{n_{2l} +n_{2j}}$ admissible.  Let
$(\beta_i)_{i=1}^t \in Ba(\ell_q)$. We wish to show that,

$$ \frac{1}{m_{2l} m_{2j}} \sum_{i=1}^t \beta_i x_i^* \in \NC.$$

\n Do this by carefully grouping the functionals.  Let
$(J_k)_{k=1}^m$ be successive intervals of integers such that
$\bigcup_{k=1}^m J_k = \{1,\ldots, t\}$, $(x_i^*)_{i\in J_k}$ is
$S_{n_{2l}}$ admissible for each $k \leq m$ and $(x^*_{\min
J_k})_{k=1}^m$ is $S_{n_{2j}}$ admissible.  Now define a sequence
$(z_k^*)_{k=1}^m$ by,

\begin{equation*}
\begin{split}
\frac{1}{m_{2l} m_{2j}} \sum_{i=1}^t \beta_i x_i^* = \frac{1}{m_{2l}
m_{2j}} \sum_{k=1}^m \sum_{i \in J_k} \beta_i x_i^* & = \frac{1}{
m_{2j}} \sum_{k=1}^m \biggl(\sum_{i \in J_k} |\beta_i|^q
\biggr)^\frac{1}{q} \frac{1}{m_{2l}} \sum_{i\in J_k}
\frac{\beta_i}{(\sum_{i \in J_k}
|\beta_i|^q )^{1/q}} x_i^* \\
 & = \frac{1}{ m_{2j}} \sum_{k=1}^m
\biggl(\sum_{i \in J_k} |\beta_i|^q \biggr)^\frac{1}{q} z_k^*
\end{split}
\end{equation*}

\n It is straightforward to check that $z_k^* \in \NC_l$ for all $k
\leq m$. The claim follows by observing that $(z_k^*)_{k=1}^m$ is
$S_{n_{2j}}$ admissible since $(x^*_{\min J_k})_{k=1}^p$ is
$S_{n_{2j}}$ admissible and $((\sum_{i\in J_k}
|\beta_i|^q)^{1/q})_{k=1}^m \in Ba(\ell_q)$.

Before proceeding further, we pause briefly to discuss the structure
of the proof of Theorem \ref{mainthm}. The proof begins by
introducing some auxiliary remarks and lemmas. Remark
\ref{decomprem} and Lemma \ref{dl} follow from the technical
definitions of the sequences $(n_i)_i$ and $(m_i)_i$ and the tree
structure of the functionals in $\NC$.  Lemma \ref{dl} is quite
specific to spaces which are $(M,N,p)$ Schreier and will be used
throughout the proof of Theorem \ref{mainthm}.   The main task at
hand is to construct a sequence of functionals in $\NC$ which are
seminormalized and satisfy an upper $\ell_p^{\omega}$ estimate with
constant 1. We do this in Lemma \ref{mainlemma}.  The construction
of these functionals is rather straightforward; it is in proving
that they possess the desired properties that we must make use of
Lemma \ref{dl} and Corollary \ref{admissiblecor}. Once we have
constructed these norming functionals (and after making a few easy
remarks) we are ready to define the operator.  This is done in a
very natural way. The fact that the operator is bounded and
non-compact follows from the the properties of the norming
functionals from which it is built.

For any functional tree $\T$ we define a function $\varphi: \T
\rightarrow \N \cup \{ 0 \}$ in the following way,

\[ \varphi(\beta) = \left\{ \begin{array}{ll}
                n_{2i} & \mbox{if $n_\beta =n_{2i}$ for some $i$} \\
                n_1  &  \mbox{if $n_\beta =n_{2i+1}$ for some $i$} \\
                0 & \mbox{if $\beta$ is terminal}
                \end{array}
        \right. \]

\begin{rem}
Let $(x_{\alpha}^*)_{\alpha \in \T}$ be a functional tree for some
$x^* \in \NC$ such that for $\alpha \in T$, $(x^*_\beta)_{\beta \in
D_{\alpha}(\T)}$ is $S_{\varphi(\T)}$ admissible. For every subset
$A$ of $\T$ consisting of pairwise incomparable nodes, the
collection $(x^*_{\alpha})_{\alpha \in A}$ is $S_d$ admissible where
$d=\max\{\sum_{\beta \prec \alpha} \varphi(\beta) : \alpha \in A\}$.
\label{decomprem}
\end{rem}

\begin{proof}

We proceed by induction on $o(\T)$.  The base step is trivial. Let
$k \geq 1$ assume the statement for $\T$ such that $o(\T) < k+1$ and
suppose $o(\T)=k+1$. Let $\alpha_0$ be the root of $\T$ and for
$\alpha \in D_{\alpha_0}(\T)$ let $\T_{\alpha}$ be the tree
corresponding to $x^*_\alpha$. For the given collection $A$ and
$\alpha \in D_{\alpha_0}(\T)$ we can define $A_\alpha =\{\beta :
\beta \in T_\alpha \cap A \}$. Notice that $A=\bigcup_{\alpha \in
D_{\alpha_0}(\T)} A_{\alpha}$ or $A=\{\alpha_0\}$. Apply the
induction hypothesis for each collection $A_\alpha$ to conclude that
$(x^*_\beta)_{\beta \in A_{\alpha}}$ is $S_{d_\alpha}$ admissible
for $d_\alpha = \max \{ \sum_{\alpha_0 \prec \gamma \prec \beta}
\varphi(\gamma) : \beta \in A_\alpha \}$. Let
$d_{\alpha_0}=\max_{\alpha \in D_{\alpha_0}(\T)}d_\alpha $.  The
block sequence $((x_\beta^*)_{\beta \in A_{\alpha}})_{\alpha \in
D_{\alpha_0}(\T)}$ is $S_{d_{\alpha_0}+\varphi(\alpha_0)}$
admissible by the convolution property of Schreier families. Finish
by observing that $d_{\alpha_0}+\varphi(\alpha_0) =
\max\{\sum_{\beta \prec \alpha} \varphi(\beta) : \alpha \in A\}$ and
$((x_\beta^*)_{\beta \in A_{\alpha}})_{\alpha \in D_{\alpha_0}(\T)}
= (x^*_\alpha)_{\alpha \in A}$.
\end{proof}

Our next lemma allows us to decompose norming functionals.
Decompositions are extremely useful when attempting to find tight
upper estimates on the norm of vectors in the space.

\begin{lem}
(Decomposition Lemma) Let $k \in \N$ and $x^* \in \NC$ such that
$\supp x^* \geq 2k$. There is an $m \in \N$, $x_1^* < \ldots < x^*_m
\in \NC$, a partition $I_1, I_2$ of $\{1,\ldots,m\}$ and scalars
$(\lambda_i)_{i=1}^m$ such that,
\begin{itemize}
\item[(a)] $x^* = \sum_{i=1}^m \lambda_i x^*_i$
\item[(b)] $x_i^* = \pm e^*_{j_i}$ for $ i \in I_1$ and $\{j_i : i \in I_1
\} \in S_{p_k-1}$.
\item[(c)] $(\sum_{i \in I_2} |\lambda_i|^q )^{1/q} \leq 2/m_{2k}$ and $(\sum_{i \in I_1 \cup I_2} |\lambda_i|^q )^{1/q} \leq 2$
\label{dl}
\end{itemize}
\end{lem}

\begin{proof}
Let $x^* \in \NC$ and $k \in \N$.  Let $\T$ be the tree
corresponding to $x^*$.  For each node $\beta$ there are
corresponding $m_{\beta}$, $n_\beta$ and $\gamma_\beta$.  Let
$\mathcal{B}$ denote the set of branches of $\T$. For each branch $b
\in \mathcal{B}$ let $\alpha(b)$ denote the node of $b$ such that
either $\alpha(b)$ is the first node $\beta$ for which
$\prod_{\alpha \prec \beta} m_{\alpha} \geq m_{2k}$ holds, or the
terminal node of $b$ if no such $\beta$ exists.  Set
$A=\{\alpha(b):b \in \mathcal{B}\}$. Notice that $A$ is a collection
of pairwise incomparable nodes intersecting every branch of
$\mathcal{B}$. Let $A_1$ denote the set of terminal nodes of $A$ and
$A_2=A \setminus A_1$. Enumerate $A$ with the set $\{1, \ldots ,
m\}$ for some $m \in \N$ and define $I_t =\{i : x_i^* \in
(x^*_{\alpha})_{\alpha \in A_t}\}$ for $t \in \{1,2\}$.  By
(\ref{decomp}) we have,

$$x^* = \sum_{\alpha \in A} \frac{\prod_{\beta \preceq \alpha}
\gamma_\beta}{\prod_{\beta \prec \alpha} m_\beta} x^*_{\alpha},
~\mbox{so set}~ \lambda_i = \frac{\prod_{\beta \preceq \alpha}
\gamma_\beta}{\prod_{\beta \prec \alpha} m_{\beta}} ~\mbox{if}~
x^*_i=x_{\alpha}^*.$$

\n It is left to verify that conditions (b) and (c) hold.  Condition
(c) follows from the fact that for each $\alpha$,
$(\gamma_{\beta})_{\beta \in D_{\alpha}(\T)} \in 2^{1/p}Ba(\ell_q)$,
and observing that,

$$\biggl( \sum_{i\in I_2} |\lambda_i|^q \biggr)^{1/q} = \biggl( \sum_{\alpha \in
A_2} \biggl|\frac{\prod_{\beta \preceq \alpha}
\gamma_\beta}{\prod_{\beta \prec \alpha} m_{\beta}}\biggr|^q
\biggr)^{\frac{1}{q}} \leq \frac{1}{m_{2k}} \biggl(\sum_{\alpha \in
A_2} \biggl|\prod_{\beta \preceq \alpha} \gamma_\beta \biggr|^q
\biggr)^\frac{1}{q} \leq \frac{2^{1/p}}{m_{2k}} <
\frac{2}{m_{2k}}.$$

\n The second part of (c) follows similarly.  The first part of (b)
follows from the definition. For the second part of (b) we employ
Remark \ref{decomprem}.  Let $\R= \bigcup_{\alpha \in A_1} \{\beta :
\beta \prec \alpha \}$. For $\alpha \in \R$ such that
$m_{\alpha}=m_{2j+1}$ for some $j \in \N$, $(x^*_{\beta})_{\beta \in
D_{\alpha}(\R)}$ is $ S_1$ and hence $S_{n_1}$ admissible. To see
this, first note that for all $\beta \in \R$, $m_\beta <m_{2k}$. By
the injectivity of the function $\sigma$, (defined in
$\NC_\infty^q$) for $\beta,\gamma \in D_\alpha(\R)$, $m_\beta \not=
m_\gamma < m_{2k}$. Since $\supp x^* \geq 2k$ we have that
$(x^*_\beta)_{\beta \in D_{\alpha}(\R)}$ is $S_1$ admissible. Thus
for $\alpha \in A_1$, $(x^*_\beta)_{\beta \in D_\alpha(\R)}$ is
$S_{\varphi(\alpha)}$ admissible. By Remark \ref{decomprem},
$(x_{\alpha}^*)_{\alpha \in A_1}$ is $S_d$ admissible where $d=
\max\{\sum_{\beta \prec \alpha} \varphi(\beta) : \alpha \in A_1\}$.

Let $\alpha \in A_1$.  We have $ \prod_{\beta \prec \alpha} m_\beta
= m_1^{b_1} \prod_{i<k} m_{2i}^{b_{2i}+5b_{2i+1}} < m_{2k}$, where
$b_j=|\{\beta : \beta \prec \alpha, m_\beta = m_j\}|$.  Apply (1.2)
for $b_1=$``$a$'' and $b_{2i}+5b_{2i+1}=$``$a_i$'', to conclude
that,

$$b_1 n_1 + \sum_{i<k} (b_{2i}+5b_{2i+1})
n_{2i} < \sum_{i<k} s_i n_{2i} =p_k.$$

\n We also have,

$$ \sum_{\beta \prec \alpha} \varphi(\beta) = \biggl( \sum_{0 \leq i<k}
b_{2i+1}\biggr) n_1 + \sum_{1 \leq i<k} b_{2i} n_{2i} < b_1n_1
+\sum_{1 \leq i<k} (b_{2i}+5b_{2i+1}) n_{2i}. $$

\n This holds for all $\alpha \in A_1$ and thus, $\max\{\sum_{\beta
\prec \alpha} \varphi(\beta) : \alpha \in A_1 \} \leq p_k -1$.
\end{proof}

\begin{cor}
Let $x^* \in \NC$ and $k \in \N$. Decompose $x^*$ as,

$$x^* = \sum_{\beta \in \max \T} \frac{\prod_{\alpha \preceq
\beta}\gamma_{\alpha}  }{\prod_{\alpha \prec \beta} m_{\alpha} }
e^*_{j_{\beta}}.$$

\n Then the set,

$$\left\{ j_{\beta} : |x^*(e_{j_\beta})| \geq \frac{ 2 \prod_{\alpha
\preceq \beta}\gamma_{\alpha}  }{m_{2k} }, j_\beta \geq 2k \right\}
$$

\n is $S_{p_k-1}$ admissible. \label{admissiblecor}
\end{cor}

\begin{proof}
For $k \in \mathbb{N}$ we can assume without loss of generality that
$\supp x^* \geq 2k$.  Apply the decomposition lemma to $x^*$ to
obtain $I_1$ and $I_2$ such that

$$x^* = \sum_{i \in I_1} \lambda_i e^*_{j_i} + \sum_{i \in I_2}
\lambda_i x^*_i $$

\n where $\{ j_i : i \in I_1 \} \in S_{p_k-1}$.  We claim that,

$$\left\{ j_{\beta} : |x^*(e_{j_\beta})| \geq \frac{ 2 \prod_{\alpha
\preceq \beta}\gamma_{\alpha} }{m_{2k} } , j_\beta \geq 2k \right\}
\subset \{ j_i : i \in I_1 \}.$$

\n If this were not the case, then for some $i_0 \in I_2$,

$$ \frac{ 2 \prod_{\alpha
\preceq \beta}\gamma_{\alpha} }{ m_{2k} } \leq |x^*(e_{j_{\beta}})|
= |\lambda_{i_0} x^*_{i_0} (e_{j_{\beta}}) | \leq |\lambda_{i_0}|.$$

\n  From the proof of the decomposition lemma,

$$ \lambda_{i_0} = \frac{ \prod_{\alpha
\preceq \beta}\gamma_{\alpha} }{\prod_{\alpha \prec \beta}
m_{\alpha} } ~\mbox{for some} ~\beta \in A_2.$$

\n For $\beta\in A_2$ we have that $\prod_{\alpha \prec \beta}
m_{\alpha} \geq m_{2k}$.  Serving as our contradiction.
\end{proof}

Before passing to the main lemma of the paper we state the following
fact concerning the existence of a particular sequence of scalars.
These scalars are called repeated hierarchy averages and were first
studied in by Argyros, Mercourakis and Tsarpalias in \cite{AMT}. These averages are defined in \cite{AB} for
$q=2$. In \cite{Ga1} a similar fact is established for $q=1$.

\begin{fact}
For any $1\leq q<\infty$ and $\vp >0$, there exist successive
subsets of $\N$, $(F_k)_{k=1}^\infty$, and scalars $(a_{k,i})_{i\in
F_k}$, such that for each $k \in \N$, $F_k \geq 2k$, $F_k \in
S_{p_k}$, $\|(a_{k,i})_{i \in F_k}\|_{q}=1$ and $(\sum_{i \in G}
|a_{k,i}|^q)^{1/q} < \vp$ for $G\in S_{p_k-1}$. \label{fact}
\end{fact}

The next lemma establishes the existence of a seminormalized block
sequence satisfying an upper $\ell_q^\omega$ estimate with constant
1 in $X_\NC$. These blocks are constructed using Fact \ref{fact} and
used to construct the desired operator on $X_\NC$.

\begin{lem}
Let $(F_k)_{k=1}^\infty$  be successive subsets of $\N$ and scalars
$(a_{k,i})_{i \in F_k}$ be such that $F_k \geq 2k$, $F_k \in
S_{p_k}$, $\|(a_{k,i})_{i \in F_k}\|_{q}=1$ and $(\sum_{i \in G}
|a_{k,i}|^q)^{1/p} < 1/m_{2k}$ for all $G \in S_{p_k-1}$  and each
$k \in \N$. The sequence of functionals $(x_k^*)_{k=1}^\infty \in
\NC$ defined by, $x_k^*=1/m_{2k} \sum_{i \in F_k} a_{k,i} e^*_i$,
are seminormalized and satisfy an upper $\ell_q^\omega$-estimate
with constant 1. \label{mainlemma}
\end{lem}

\begin{proof}
We start by making an observation concerning the decomposition of
each $x^*_k$. For fixed $k$, and $k_0 \leq k$ write $F_k
=\bigcup_{r=1}^{d_k} J_{k,r}$ such that $J_{k,1} < \ldots <
J_{k,d_k}$, each $J_{k,r}$ is $S_{p_k-p_{k_0}}$ admissible and
$(J_{k,r})_{r=1}^{d_k}$ is $S_{p_{k_0}}$ admissible (we can do this
because $F_k$ is $S_{p_k}$ admissible). Then,

$$x^*_k=\frac{1}{m_{2k_0}} \sum_{r=1}^{d_k} \biggl(\sum_{i \in J_{k,r}}
|a_{k,i}|^q \biggr)^\frac{1}{q} z^*_{k,r} ~\mbox{for}~
z^*_{k,r}=\frac{m_{2k_0}}{m_{2k}} \sum_{i \in J_{k,r}}
\frac{a_{k,i}}{(\sum_{i \in J_{k,r}} |a_{k,i}|^q )^{1/q}} e^*_i.$$

\n Since $m_{2k_0} / m_{2k} = 1/ \prod_{2k_0 \leq \ell <2k}
m_\ell^{s_\ell}$ and $(e^*_i)_{i \in J_{k,r}}$ is $S_{p_k-p_{k_0}}$
admissible, we conclude by (1.3) that $z^*_{k,r} \in \NC$ for all $r
\leq d_k$. Since, $\min J_{k,r} =\min \supp z^*_{k,r}$, we have that
$(z^*_{k,r})_{r=1}^{d_k}$ is $S_{p_{k_0}}$ admissible.

We now show that $(x_k^*)_k$ satisfies and upper $\ell_q^\omega$
estimate with constant 1.  For starters, let $k_0 \in \N$ and $F
\subset \N$ with $F \geq k_0$ such that $(x_k^*)_{k \in F}$ is
$S_{f_{k_0}}$ admissible. For every $k \in F$ we apply the above
(since $F \geq k_0$) to define $(z^*_{k,r})_{r=1}^{d_k}$.  The block
sequence $((z^*_{k,r})_{r=1}^{d_k})_{k \in F}$ is
$S_{p_{k_0}+f_{k_0}}$ admissible, by the convolution property of
Schreier families. Hence it is $S_{n_{2k_0}}$ admissible by (1.1)
and the hereditary property of Schreier families.  To conclude, it
suffices to let $(\beta_k)_{k \in F} \in Ba(\ell_q)$ and show that
$\sum_{i\in F} \beta_i x^*_i \in \NC$. We do this by observing the
following equality,

\begin{equation*}
\begin{split}
\sum_{k \in F} \beta_k x^*_k & = \sum_{k \in F} \beta_k
\frac{1}{m_{2k_0}} \sum_{r=1}^{d_k} \biggl( \sum_{i \in J_{k,r}}
|a_{k,i}|^q \biggr)^\frac{1}{q} z^*_{k,r} \\
    & = \frac{1}{m_{2k_0}} \sum_{k\in F} \sum_{r=1}^{d_k} \beta_k \biggl(
    \sum_{i \in J_{k,r}} |a_{k,i}|^q \biggr)^\frac{1}{q} z^*_{k,r}.
\end{split}
\end{equation*}

\n Since $(\beta_k ( \sum_{i \in J_{k,r}} |a_{k,i}|^q )^\frac{1}{q}
)_k \in Ba(\ell_q)$ and  $((z^*_{k,r})_{r=1}^{d_k})_{k \in F}$ is
$S_{n_{2k_0}}$ admissible, it follows that $ \sum_{k\in F} \beta_k
x^*_k \in \NC$.  Thus, $(x_k^*)_k$ satisfies a upper
$\ell_q^\omega$-estimate with constant 1.

To show that $(x_k^*)_k$ is seminormalized, it suffices to find a
uniform lower bound.  For each $k$, define $x_k=\sum_{j \in F_k}
a_{k,j}^{q/p} e_j$. It suffices to show that $\|x_k\| \leq
26/m_{2k}$. From this it follows easily that $\|x_k^*\| \geq 1/26$
for all $k \in \N$. Let $x^* \in \NC$ be an arbitrary norming
functional which we may assume without loss of generality satisfies
$\supp x^* \geq 2k$ (since $F_k \geq 2k$). By applying the
decomposition lemma for $k \in \N$ and $x^*$, we can estimate
$\|x_k\|$ from above as follows,

\begin{equation*}
\begin{split}
|x^*(x_k)| & \leq \biggl| \sum_{i \in I_1} \lambda_i e^*_{j_i} (x_k)
\biggr| + \biggl| \sum_{i \in I_2} \lambda_i y^*_i (x_k) \biggr| \\
    & \leq \sum_{i \in I_1} |\lambda_i| |a_{k,j_i}|^\frac{q}{p} +
    \sum_{i \in I_2} |\lambda_i| 12 \biggl( \sum_{\{j : j \in \supp
    y_i \cap F_k \}} |a_{k,j}|^q \biggr)^\frac{1}{p} \\
    & \leq \biggl( \sum_{i \in I_1} |\lambda_i|^q
    \biggr)^\frac{1}{q} \biggl( \sum_{i \in I_1} |a_{k,j_i}|^q
    \biggr)^\frac{1}{p} + 12 \biggl( \sum_{i \in I_2} |\lambda_i|^q
    \biggr)^\frac{1}{q} \biggl( \sum_{i \in I_2}  \sum_{\{j : j \in \supp
    y_i \cap F_k \}} |a_{k,j}|^q \biggr)^\frac{1}{p} \\
    & \leq 2\frac{1}{m_{2k}} + 12\frac{2}{m_{2k}} = \frac{26}{m_{2k}}.
\end{split}
\end{equation*}

\n The first inequality follows from the decomposition lemma and the
triangle inequality. The second inequality follows from the triangle
inequality, the definition of $x_k$, and Remark \ref{ulp}. The third
follows from two applications of H\"{o}lders inequality. For the
last inequality we used condition (c) of the decomposition lemma,
the fact that $(j_i)_{i \in I_1}$ is $S_{p_k-1}$ admissible (by
condition (b) of the decomposition lemma) and the definition of
$(a_{k,i})_{i \in F_k}$.  This concludes the proof.
\end{proof}

We make two final remarks before proceeding with the proof of the
main theorem.

\begin{rem}
Let $(y_i^*)_i$ be the even subsequence of the seminormalized block
sequence $(x^*_i)_i$ satisfying an upper $\ell_q^\omega$ estimate
with constant $1$ defined in Lemma \ref{mainlemma}. Let $k \in \N$,
$F \subset \N$, with $F\geq k$ such that $(y^*_i)_{i\in F}$ is
$S_{n_{2k}}$ admissible.  Then $\|\sum_{i \in F} \beta_i y^*_i \|
\leq 1$ for all $(\beta_i)_i \in Ba(\ell_q)$. \label{sub}
\end{rem}

\begin{proof}
Let $k \in \N$, $F$ be a subset of $\N$, with $F\geq k$ and
$(y^*_i)_{i\in F}$ being $S_{n_{2k}}$ admissible.  Set $G=\{i :
i=2j, ~ j \in F\}$ and note that $(y_i)_{i \in F} = (x^*_i)_{i \in
G}$. Since $F \geq k$ we have, $i \geq k+1 $ for all $i \in G$.
Since $(x^*_i)_i$ satisfies an upper $\ell_q^\omega$, estimate $G
\geq k+1$ and $(x_i^*)_{i \in G}$ is $S_{n_{2k}}$ admissible an thus
$S_{f_{k+1}}$ admissible, we have,

$$\biggl\|\sum_{i \in F} \beta_i y_i^* \biggr\| = \biggl\| \sum_{ i \in G} \beta_i
x^*_i \biggr\| \leq 1 .$$

\n This concludes the proof.
\end{proof}

\begin{rem}
Let $(y_i^*)_i$ be the subsequence from Remark \ref{sub}. For every
$x \in S(X)$, $k \in \N$, $F \subset \N$ with $F \geq k$ and
$(y^*_i)_{i \in F}$ being $S_{n_{2k}}$ admissible we have
$(y_i^*(x))_{i\in F} \in Ba(\ell_p)$. \label{firstrem}
\end{rem}

\begin{proof}
\n Let $x \in S(X)$,  $k \in \N$ and $F \subset \N$ with $F \geq k$
such that $(y^*_i)_{i \in F}$ is $S_{n_{2k}}$ admissible. By Remark
\ref{sub} for all $(\beta_i)_{i \in F} \in Ba(\ell_q)$ we have $\|
\sum_{i \in F} \beta_i y^*_i \| \leq 1$.   Apply this for,

$$\beta_i = \frac{|y_i^*(x)|^{p/q}\sign(y_i^*(x))}{(\sum_{j\in F}
|y^*_j(x)|^p)^{1/q}}$$

\n and estimate $\|\sum_{i \in F} \beta_i y^*_i \|$ from below with
$x$.
\end{proof}

\n ({\em Proof of Theorem \ref{mainthm}}) We are now ready to define
the desired operator on $X_\NC$. Let $(y_i^*)_i$ be the
seminormalized block sequence from Remark \ref{sub}. For $x \in
c_{00}$ define the operator $T: c_{00} \rightarrow c_{00}$ by $T x =
\sum_{i=1}^\infty y_i^*(x)e_i$.  Once we show that $T$ is a bounded
operator it can be extended as on operator defined on $X_\NC$.

Since $(y_i^*)_i$ is a seminormalized block sequence it follows that
$T$ is non-compact.  In the case that $X_\NC$ is an HI space, $T$
must be strictly singular. Since $dim(\Ker T)=\infty$, if there was
an infinite dimensional subspace $Y$ of $X_{\NC}$ such that $T|_{Y}$
was an isomorphism. $Y+\Ker(T)$ would be a direct sum. Contradicting
the fact that $X_\NC$ is HI. (It is known that the spaces
constructed in \cite{DM} have few operators. Using similar
techniques it can be further shown that the space constructed in
\cite{AB} has few operators.)

Our final task is to demonstrate that $T$ is bounded.  Let $x \in
S(X)$ and $x^* \in \NC$. If $x^*=\pm e^*_j$ for some $j$ then
$|x^*(Tx)| \leq 1$. Thus assume $x^* \in \NC$ such that $|\supp x^*|
>1$.  Suppose $x^*$ has the following decomposition,

$$x^* = \sum_{\beta \in \max \T } \frac{\prod_{\alpha \preceq
\beta}\gamma_{\alpha} }{\prod_{\alpha \prec \beta} m_{\alpha} }
e^*_{j_{\beta}}.$$

\n Define,

$$H_2 = \biggl\{ j_{\beta} : \frac{2 \prod_{\alpha \preceq \beta}
\gamma_\alpha}{m_{4}} \leq |x^*(e_{j_\beta})| < \frac{\prod_{\alpha
\preceq \beta} \gamma_\alpha}{m_1}\biggr\},$$

\n For $k > 2$ define,

$$H_k = \biggl\{ j_{\beta} : \frac{2 \prod_{\alpha \preceq \beta}
\gamma_\alpha}{m_{2k}} \leq |x^*(e_{j_\beta})| < \frac{2
\prod_{\alpha \preceq \beta} \gamma_\alpha}{m_{2(k-1)}}\biggr\}.$$

\n For $k > 2$ define, $G_k = \{ j_{\beta} \in H_k : j_{\beta} \geq
2k \}. $  Clearly, $\supp x^* =\bigcup_{k=2}^\infty H_k$. Apply
Corollary \ref{admissiblecor} to deduce that $G_k \in S_{p_k-1}$. By
(1.1) and (ii), $ G_k \in S_{n_{2k}}$. By the spreading property of
Schreier families $(y_i)_{i \in G_k}$ is $S_{n_{2k}}$ admissible for
all $k$. For each $k$, apply Remark \ref{firstrem} to deduce that,

\begin{equation}
 \biggl( \sum_{i \in G_k} |y_i^*(x)|^p \biggr)^\frac{1}{p} \leq 1.
 \quad\quad \biggl( \sum_{i \in H_2} |y_i^*(x)|^p \biggr)^\frac{1}{p} \leq 1  \label{last}
 \end{equation}

\n Estimate $|x^*(Tx)|$ from above in the following way,

\begin{equation*}
\begin{split}
x^*\biggl(\sum_{i=1}^\infty y_i^*(x)e_i \biggr) & \leq \sum_{i\in
H_2} |y^*_i(x)||x^*(e_i)| + \sum_{k=3}^\infty \biggl[\sum_{i\in G_k}
|y^*_i(x)||x^*(e_i)| + \sum_{i\in
H_k\setminus G_k} |y_i^*(x)||x^*(e_i)| \biggr] \\
  & \leq \biggl(\sum_{i \in H_2} |y_i^*(x)|^p \biggr)^\frac{1}{p} \biggl(
  \sum_{i\in H_2} |x^*(e_i)|^q \biggr)^\frac{1}{q} + \sum_{k=3}^\infty \biggl[ \biggl(\sum_{i \in G_k} |y_i^*(x)|^p
    \biggr)^\frac{1}{p} \biggl(
  \sum_{i\in G_k} |x^*(e_i)|^q \biggr)^\frac{1}{q} \\
  &\quad + \biggl(\sum_{i \in H_k\setminus G_k} |y_i^*(x)|^p \biggr)^\frac{1}{p} \biggl(
  \sum_{i\in H_k\setminus G_k} |x^*(e_i)|^q \biggr)^\frac{1}{q}
  \biggr] \\
  & < \biggl( \sum_{j_\beta \in H_2} \biggl| \prod_{\alpha \preceq
  \beta} \gamma_\alpha \biggr|^q \biggr)^\frac{1}{q} \frac{1}{m_1} + \sum_{k=3}^\infty \biggl[ \biggl(
  \sum_{j_\beta \in G_k} \biggl| \prod_{\alpha \preceq \beta} \gamma_\alpha \biggr|^q
  \biggr)^\frac{1}{q} \frac{2}{m_{2(k-1)}}  \\
  & \quad + \biggl(
  \sum_{j_\beta \in H_k \setminus G_k} \biggl| \prod_{\alpha \preceq \beta} \gamma_\alpha \biggr|^q
  \biggr)^\frac{1}{q}\frac{2(k-1)}{m_{2(k-1)}} \biggr] \\
  & \leq \frac{2}{m_1} + \sum_{k=3}^\infty \frac{4}{m_{2(k-1)}} +
  \frac{4(k-1)}{m_{2(k-1)}} = M.
\end{split}
\end{equation*}

\n The first inequality follows from the triangle inequality and the
definitions of $H_k, G_k$. For the second, we apply H\"{o}lders
inequality to each of the terms. For the first and second terms of
the third inequality we used (\ref{last}) and the definition of
$H_k$. For the third term of the third inequality  we used the fact
that $|x^*(e_i)|\leq 1$ for all $i$, $|H_k \setminus G_k| \leq k-1$
and the definition of $H_k$. For the final inequality we used the
fact that $(|\prod_{\alpha \preceq \beta}\gamma_\alpha |)_{\beta \in
A} \in 2Ba(\ell_q)$ for $A= H_k, G_k$ or $H_k \setminus G_k$. Thus,
$\|T\| \leq \max\{M,1\}$.

We conclude by noting that $\ell_\infty$ embeds isomorphically into
$\mathcal{L}(X_\NC)$ via the mapping

$$(a_i)_{i=1}^\infty \longmapsto SOT - \lim_{n \rightarrow \infty} \sum_{i=1}^n a_i y^*_i \otimes e_i. $$

\n Here ``$SOT-\lim$'' denotes the strong operator topology limit.
To see that this is a bounded isomorphism one merely follows, almost
identically, the previous calculation.

{\footnotesize
\noindent
Department of Mathematics and Computer Science, Amherst College, Amherst, MA 01002, \\
\noindent kbeanland@amherst.edu }

\end{document}